\documentclass[12pt,leqno]{article}
\usepackage{amsfonts}
\usepackage{listings}
\usepackage{fullpage}
\usepackage{amssymb}
\usepackage{amsmath}
\usepackage{amsthm}
\newtheorem{theorem}{Theorem}[section]

\newtheorem{corollary}{Corollary}[section]

\newtheorem{lemma}{Lemma}[section]

\numberwithin{equation}{section}
\newtheorem*{theorema}{Theorem A}
\newtheorem*{theoremb}{Theorem B}

\begin{document}

\title{\large\bf On least periods of a class of arithmetic functions    \footnote{This work is supported by the National Natural Science
Foundation of China (Grants No. 12371005 and 12671008).}}
\date{}
\author{Wen-Jie Huang$^1$ and Quan-Hui Yang$^{1,2}$
\footnote{Emails:~2730699552@qq.com, yangquanhui01@163.com.}\\
\small $^1$School of Mathematical Sciences\\
\small $^2$Ministry of Education Key Laboratory for NSLSCS\\
\small Nanjing Normal University, Nanjing 210023, P. R. CHINA}
\maketitle

\begin{abstract} Let $k, a, c$ be positive integers, let $b$ be a nonnegative integer, and let $f$ be a nonvanishing complex-valued multiplicative function. Qian, Tan and Hong proved that
$$
G_f(n)=\frac{\prod\limits_{i=0}^k f(b+a(n+i c))}{f\left(\operatorname*{lcm}\limits_{0 \leq i \leq k}\{b+a(n+i c)\}\right)}
$$
is periodic and asked for its least positive period. In this paper, we answer this problem by giving an explicit formula for this period in terms of finitely many values of $f$ at prime powers. In particular, we determine the least periods for $f(n)=\sigma_\alpha(n)=\sum\limits_{d \mid n} d^\alpha$ for every real $\alpha$, and for $f(n)=n^{\varepsilon}$ for every real $\varepsilon$.



\end{abstract}

\section{Introduction}

Let $\mathbb{N}=\{0,1,2, \ldots\}$ and $\mathbb{Z}^{+}=\{1,2, \ldots\}$, and let $\mathbb{Z}$ and $\mathbb{C}$ denote the integers and complex numbers. We use gcd and lcm for the greatest common divisor and least common multiple, $p$ for a prime, and $\# S$ for the number of elements of a finite set $S$. For $m \in$ $\mathbb{Z} \backslash\{0\}, v_p(m)$ denotes the exponent of $p$ in $|m|$, and we set $v_p(0)=+\infty$. For $m \neq 0$, we write $p^{v_p(m)} \| m$ to mean that $p^{v_p(m)} \mid m$ and $p^{v_p(m)+1} \nmid m$. For $k \in \mathbb{N}$ and $n \in \mathbb{Z}^{+}$, Farhi~\cite{Farhi} introduced the arithmetic function

$$
g_k(n)=\frac{n(n+1) \cdots(n+k)}{\operatorname{lcm}(n, n+1, \ldots, n+k)}
$$
and proved that $k!$ is a period. Hong and Yang~\cite{HongYang} showed that $L_k$ is a period, where $L_0=1$ and $L_k=\operatorname{lcm}(1, \ldots, k)$ for $k \geq 1$. Farhi and Kane determined the least positive period of $g_k$.

\begin{theorema}[{\cite{FarhiKane}, Theorem~3.2}] For every $k \in \mathbb{N}$, the least positive period $P_k$ of $g_k$ is given by

$$
P_k=\prod\limits_{\substack{p \text { prime } \\ p \leq k}} p^{\delta_p(k)},
$$
where

$$
\delta_p(k)= \begin{cases}0, & v_p(k+1) \geq \max\limits_{1 \leq i \leq k} v_p(i), \\ \max\limits_{1 \leq i \leq k} v_p(i), & \text { otherwise } .\end{cases}
$$
\end{theorema}

Here and below, an empty product has value 1.
For integers $k \geq 0, a \geq 1$ and $b \geq 0$, put

$$
g_{k, a, b}(n)=\frac{(b+n a)(b+(n+1) a) \cdots(b+(n+k) a)}{\operatorname{lcm}(b+n a, b+(n+1) a, \ldots, b+(n+k) a)} \quad\left(n \in \mathbb{Z}^{+}\right)
$$
and $e_{p, k}=\left\lfloor\log _p k\right\rfloor=\max\limits_{1 \leq i \leq k} v_p(i)$ for $k \geq 1$. Hong and Qian extended the result of Farhi and Kane to arithmetic progressions.

\begin{theoremb}[{\cite{HongQian}, Theorem~1.2}] The function $g_{k, a, b}$ is periodic. If $\operatorname{gcd}(a, b)=1$, its least positive period is $Q_{k, a}$, where

$$
Q_{k, a}:=\frac{L_k}{\delta_{k, a} \prod\limits_{\substack{q \text { prime } \\ q \mid \operatorname{gcd}\left(a, L_k\right)}} q^{e_{q, k}}},
$$
and

$$
\delta_{k, a}:= \begin{cases}p^{e_{p, k}}, & \text { if } p \nmid a \text { and } v_p(k+1) \geq e_{p, k} \text { for some prime } p \leq k, \\ 1, & \text { otherwise } .\end{cases}
$$
If $\operatorname{gcd}(a, b)>1$, its least positive period is $Q_{k, a^{\prime}}$, where
$a^{\prime}=\frac{a}{\operatorname{gcd}(a, b)}$.
\end{theoremb}

The prime in the first case of the definition of $\delta_{k, a}$ is unique whenever it exists, by Farhi and Kane~\cite[Proposition~3.3]{FarhiKane}.

For positive integers $h$ and $k$, Ji and Ji~\cite[Theorem~1.7]{JiJi} proved that the function

$$
\frac{n(n+h) \cdots(n+k h)}{\operatorname{lcm}(n, n+h, \ldots, n+k h)} \quad\left(n \in \mathbb{Z}^{+}\right)
$$
has least positive period $h P_k$. They also obtained the least positive period $Q_{k, a}$ of $g_{k, a, b}$ when $k \geq 1$ and $\operatorname{gcd}(a, b)=1$~\cite[Corollary~1.8]{JiJi}.

Let $k, a, c \in \mathbb{Z}^{+}, b \in \mathbb{N}$, and let $f: \mathbb{Z}^{+} \rightarrow \mathbb{C} \backslash\{0\}$ be multiplicative. Define

\begin{equation}\label{eq:G}G_f(n)=\frac{\prod\limits_{i=0}^k f(b+a(n+i c))}{f\left(\operatorname*{lcm}\limits_{0 \leq i \leq k}\{b+a(n+i c)\}\right)} \quad\left(n \in \mathbb{Z}^{+}\right) \end{equation}

Qian, Tan and Hong~\cite[Problem~1.1]{QianTanHong} asked whether $G_f$ is periodic and, if so, what its least positive period is. They proved that $c L_k$ is a period and determined the least period when $f$ is Euler's totient function. They also asked for the cases $f(n)=\sigma_\alpha(n)=\sum\limits_{d \mid n} d^\alpha$ with $\alpha \in \mathbb{N}$ and $f(n)=n^{\varepsilon}$ with $\varepsilon \in \mathbb{R}$.

In this paper, we determine the least positive period of $G_f$ for every nonvanishing complex-valued multiplicative function $f$. The formula uses only finitely many values of $f$ at prime powers. It also gives the least period for $\sigma_\alpha$ for every real $\alpha$, as well as for $n^{\varepsilon}$ for every real $\varepsilon$.

Put
$$
d=\operatorname{gcd}(a, b), \quad A=\frac{a}{d}, \quad B=\frac{b}{d}, \quad L=L_k=\operatorname{lcm}(1, \ldots, k) .
$$
Then $\operatorname{gcd}(A, B)=1$. For each prime $p, v_p(L)=\left\lfloor\log _p k\right\rfloor$, so $p^{v_p(L)} \leq k<p^{v_p(L)+1}$. For $p \mid c L$ with $p \nmid A$, define the ratios of adjacent prime power values by

$$
R_{p, j}=\frac{f\left(p^{v_p(d)+v_p(c)+j}\right)}{f\left(p^{v_p(d)+v_p(c)+j-1}\right)} \quad(j \geq 1) .
$$

We now define
$$
M_p=\max \left(\left\{j: v_p(k+1)<j \leq v_p(L), R_{p, j} \neq 1\right\} \cup\{0\}\right) .
$$
Thus $M_p$ is the largest such index greater than $v_p(k+1)$ with $R_{p, j} \neq 1$, and $M_p=0$ if no such index exists. When $M_p=0$, we write

$$
K_p=f\left(p^{v_p(d)+v_p(c)}\right)^k \prod\limits_{j=1}^{\min \left(v_p(L), v_p(k+1)\right)} R_{p, j}^{\frac{k+1}{p^j}-1} .
$$
All exponents in this product are nonnegative integers.

\begin{theorem}\label{thm:main} Let $G_f$ be defined by \eqref{eq:G}. For each prime $p \mid c L$ with $p \nmid A$, put

$$
e_p= \begin{cases}v_p(c)+M_p, & M_p>0, \\ \max \left(\left\{j+1: 0 \leq j<v_p(c), f\left(p^{v_p(d)+j}\right)^k \neq K_p\right\} \cup\{0\}\right), & M_p=0 .\end{cases}
$$
Then the least positive period of $G_f$ is

$$
P_f=\prod\limits_{\substack{p \mid c L \\ p \nmid A}} p^{e_p} .
$$
\end{theorem}

Here $e_p$ is the exponent of $p$ in the least period. For each relevant prime $p$, the formula uses only $f\left(p^{v_p(d)+j}\right)$ with $0 \leq j \leq v_p(c)+v_p(L)$. It determines the least period without computing the values of $G_f$ over a complete period.

Theorem~\ref{thm:main} recovers the result of Hong and Qian by taking $f(m)=m$ and $c=1$. Indeed, $G_f=g_{k, a, b}$, and for each $p \mid L$ with $p \nmid A$, we have $R_{p, j}=p$ and

$$
e_p=M_p= \begin{cases}v_p(L), & v_p(L)>v_p(k+1), \\ 0, & v_p(L) \leq v_p(k+1) .\end{cases}
$$
Since $e_{p, k}=v_p(L)$, it follows that

$$
P_f=\prod\limits_{\substack{p \mid L \\ p \nmid A \\ v_p(L)>v_p(k+1)}} p^{v_p(L)}=\frac{L}{\delta_{k, A} \prod\limits_{\substack{q \text { prime } \\ q \mid \operatorname{gcd}(A, L)}} q^{e_{q, k}}}=Q_{k, A} .
$$

The cases $\operatorname{gcd}(a, b)=1$ and $\operatorname{gcd}(a, b)>1$ correspond to $A=a$ and $A=a^{\prime}$, respectively. The case $k=0$ follows directly from $g_{0, a, b} \equiv 1$. Taking $a=1$ and $b=0$ also recovers the result of Farhi and Kane.

We say that a positive-valued arithmetic function $f$ is {\em strictly prime-power monotonic} if, for each fixed prime $p$, the sequence $\left(f\left(p^m\right)\right)_{m \geq 0}$ is strictly monotonic. The direction of monotonicity may depend on $p$.

\begin{corollary}\label{cor:monotone} If $f$ is strictly prime-power monotonic, then the least positive period of $G_f$ is

$$
\prod\limits_{\substack{p \mid c \\ p \nmid A}} p^{v_p(c)} \prod\limits_{\substack{p \mid L \\ p \nmid A \\ v_p(L)>v_p(k+1)}} p^{v_p(L)} .
$$
\end{corollary}

Since $\sigma_\alpha$ is strictly prime-power monotonic for every real $\alpha$, and so is $n^{\varepsilon}$ for every real $\varepsilon \neq 0$, we obtain the following corollary.

\begin{corollary}\label{cor:standard} For real $\alpha$, let $\sigma_\alpha(n)=\sum\limits_{d \mid n} d^\alpha$. Then the least positive period of $G_f$ with $f=\sigma_\alpha$ is the period given in Corollary~\ref{cor:monotone}. The same holds for $f(n)=n^{\varepsilon}$ with any real $\varepsilon \neq 0$. For $\varepsilon=0$, the least positive period is 1.
\end{corollary}

\section{Preliminary lemmas}

It is well known that the least positive period divides every positive period. Indeed, if $T_0$ is the least period of an arithmetic function $H(n)$ and $T=q T_0+r$ is any period, with $0 \leq r<T_0$, then $H(n+r)=H\left(n+r+q T_0\right)=H(n)$, so $r=0$.

\begin{lemma}\label{lem:product-period} Let $q_1, q_2, \ldots, q_m$ be pairwise coprime positive integers. For each $1 \leq i \leq m$, let $H_i: \mathbb{Z}^{+} \rightarrow \mathbb{C} \backslash\{0\}$ have least positive period $q_i$. Then the function $H(n)=\prod\limits_{i=1}^m H_i(n)$ has least positive period $\prod\limits_{i=1}^m q_i$.
\end{lemma}

\begin{proof} Put $T=\prod\limits_{i=1}^m q_i$. Clearly, $T$ is a period of $H(n)$. Conversely, let $T^{\prime}$ be any positive period of $H(n)$. For each $i$, put $T_i=\frac{T}{q_i}$. Since $T^{\prime} T_i$ is also a period of $H(n)$ and a multiple of $q_j$ for every $j \neq i$, cancelling the nonzero factors

$$
\prod\limits_{\substack{j=1 \\ j \neq i}}^m H_j\left(n+T^{\prime} T_i\right)=\prod\limits_{\substack{j=1 \\ j \neq i}}^m H_j(n)
$$
from $H\left(n+T^{\prime} T_i\right)=H(n)$ gives $H_i\left(n+T^{\prime} T_i\right)=H_i(n)$ for any positive integer $n$. Thus $q_i \mid T^{\prime} T_i$. Since $\operatorname{gcd}\left(q_i, T_i\right)=1$, we obtain $q_i \mid T^{\prime}$ for every $i$, and so $T \mid T^{\prime}$.

Therefore $T$ is the least positive period of $H(n)$.
\end{proof}

\begin{lemma}\label{lem:counting} Let $p$ be a prime, let $u$ be an integer with $\operatorname{gcd}(u, p)=1$, and let $k, j$ be positive integers such that $p^j \leq k$. For each integer $z$, define

$$
N_j(z)=\#\left\{i \in \mathbb{Z}: 0 \leq i \leq k, p^j \mid z+u i\right\} .
$$

If $p^j \mid k+1$, then $N_j(z)=\frac{k+1}{p^j}$ for every integer $z$. Otherwise, $N_j(z)$ takes two consecutive integer values and has least positive period $p^j$.
\end{lemma}

\begin{proof} Write $k+1=\ell p^j+r$, where $0 \leq r<p^j$. Since $\operatorname{gcd}\left(u, p^j\right)=1$, each block of $p^j$ consecutive indices $i$ contains exactly one solution of $z+u i \equiv 0\left(\bmod p^j\right)$. Thus $N_j(z)=\ell=\frac{k+1}{p^j}$ if $r=0$, and $N_j(z) \in\{\ell, \ell+1\}$ otherwise.

Next we suppose that $r>0$. Let $\mathbf{1}_{q \mid t}$ be 1 if $q \mid t$, and 0 otherwise. Then
$$
D(z):=N_j(z+u)-N_j(z)=\mathbf{1}_{p^j \mid z+u(k+1)}-\mathbf{1}_{p^j \mid z} .
$$

Since $p^j \nmid u(k+1)$, we have $D(z)=-1$ if and only if $p^j \mid z$. If $T$ is a positive period of $N_j(z)$, then
$$
D(T)=N_j(T+u)-N_j(T)=N_j(u)-N_j(0)=D(0)=-1,
$$
so $p^j \mid T$. Since $p^j$ is also a period, it follows that $p^j$ is the least positive period. Since $p^j>1$, the function $N_j(z)$ is nonconstant, so both values $\ell$ and $\ell+1$ occur.
\end{proof}

\section{Proofs}

\begin{proof}[Proof of Theorem~\ref{thm:main}] By $b+a(n+i c)=d[B+A(n+i c)]$, we have

$$
\begin{aligned}
b+a(n+i c) & =\prod\limits_p p^{v_p(d)+v_p(B+A(n+i c))}, \\
\operatorname*{lcm}\limits_{0 \leq i \leq k}\{b+a(n+i c)\} & =\prod\limits_p p^{v_p(d)+\max\limits_{0 \leq i \leq k} v_p(B+A(n+i c))}
\end{aligned} .
$$
Thus, by the multiplicativity of $f$, we obtain $G_f(n)=\prod\limits_p F_p(n)$, where

\begin{equation}\label{eq:prime-factor}
F_p(n)=\frac{\prod\limits_{i=0}^k f\left(p^{v_p(d)+v_p(B+A(n+i c))}\right)}{f\left(p^{v_p(d)+\max\limits_{0 \leq i \leq k} v_p(B+A(n+i c))}\right)} .
\end{equation}

For each fixed $n$, only finitely many factors differ from 1. Thus, although the product above is written over all primes, it is in fact a finite product, and there is no need to consider convergence.

Next, we consider the primes $p$ in the product according to the following cases.

{\bf Case 1.} $p \mid A$. Since $\operatorname{gcd}(A, B)=1$, we have $p \nmid B$. Hence $v_p(B+A(n+i c))=0$ for every $i$, and $F_p(n)=f\left(p^{v_p(d)}\right)^k$.

{\bf Case 2.} $p \nmid A$ and $p \nmid c L$. Then $p>k$ and $p \nmid c$. If two terms $B+A(n+i c)$ and $B+A(n+j c)$ are divisible by $p$ for some $0 \leq i<j \leq k$, then their difference $A c(j-i)$ is also divisible by $p$. Since $p \nmid A c$ and $0<j-i \leq k<p$, this is a contradiction.

Thus at most one term $B+A(n+i c)$ is divisible by $p$. If there is such a term, its corresponding factor cancels with the denominator in \eqref{eq:prime-factor}.

Therefore, in both Cases 1 and 2, we have $F_p(n)=f\left(p^{v_p(d)}\right)^k$, and so

\begin{equation}\label{eq:finite-product}G_f(n)=\left(\prod\limits_p f\left(p^{v_p(d)}\right)^k\right) \prod\limits_{\substack{p \mid c L \\ p \nmid A}} \frac{F_p(n)}{f\left(p^{v_p(d)}\right)^k}=f(d)^k \prod\limits_{\substack{p \mid c L \\ p \nmid A}} \frac{F_p(n)}{f\left(p^{v_p(d)}\right)^k} .
\end{equation}

It remains to determine the least period of the product over the fixed finite set of primes $p \mid c L$ with $p \nmid A$. Fix such a prime $p$. Since $v_p(A c)=v_p(c)$, write $A c=p^{v_p(c)} u$ with $u \in \mathbb{Z}^{+}$and $p \nmid u$. We first derive an expression for $F_p(n)$.

If $0 \leq v_p(B+A n)<v_p(c)$, then

$$
v_p(B+A(n+i c))=v_p\left(B+A n+i p^{v_p(c)} u\right)=v_p(B+A n)
$$
for every $i$, and hence $F_p(n)=f\left(p^{v_p(d)+v_p(B+A n)}\right)^k$.
Now we suppose that $v_p(B+A n) \geq v_p(c)$. Write $B+A n=p^{v_p(c)} z$ with $z \in \mathbb{Z}^{+}$. Since
$$
v_p\left(B+A n+i p^{v_p(c)} u\right)=v_p(c)+v_p(z+i u),
$$
we arrange
$$
v_p(z), v_p(z+u), \ldots, v_p(z+k u)
$$
in nondecreasing order as $w_0 \leq \cdots \leq w_k$, including equal values. As in Lemma~\ref{lem:counting}, let
$$
N_j(z)=\#\left\{i: 0 \leq i \leq k, p^j \mid z+i u\right\} \quad(j \geq 1) .
$$

For $1 \leq j \leq w_{k-1}$, exactly $N_j(z)-1$ of $w_0, \ldots, w_{k-1}$ are at least $j$. Then

$$
\begin{aligned}
F_p(n) & =\frac{\prod\limits_{i=0}^k f\left(p^{v_p(d)+v_p(B+A(n+i c))}\right)}{f\left(p^{v_p(d)+\max\limits_{0 \leq i \leq k} v_p(B+A(n+i c))}\right)}=\frac{\left(\prod\limits_{i=0}^{k-1} f\left(p^{v_p(d)+v_p(c)+w_i}\right)\right) f\left(p^{v_p(d)+v_p(c)+w_k}\right)}{f\left(p^{v_p(d)+v_p(c)+w_k}\right)} \\
& =\prod\limits_{i=0}^{k-1} f\left(p^{v_p(d)+v_p(c)+w_i}\right)=\prod\limits_{i=0}^{k-1}\left(f\left(p^{v_p(d)+v_p(c)}\right) \prod\limits_{j=1}^{w_i} \frac{f\left(p^{v_p(d)+v_p(c)+j}\right)}{f\left(p^{v_p(d)+v_p(c)+j-1}\right)}\right) \\
& =f\left(p^{v_p(d)+v_p(c)}\right)^k \prod\limits_{i=0}^{k-1} \prod\limits_{j=1}^{w_i} R_{p, j}=f\left(p^{v_p(d)+v_p(c)}\right)^k \prod\limits_{j=1}^{w_{k-1}} \prod\limits_{\substack{i=0 \\
w_i \geq j}}^{k-1} R_{p, j} \\
& =f\left(p^{v_p(d)+v_p(c)}\right)^k \prod\limits_{j=1}^{w_{k-1}} R_{p, j}^{N_j(z)-1} .
\end{aligned}
$$

Next we extend the upper limit of the product above from $w_{k-1}$ to $v_p(L)$. We shall prove $w_{k-1} \leq v_p(L)$ and $N_j(z)=1$ for $w_{k-1}<j \leq v_p(L)$.

Since $p^{v_p(L)} \leq k<p^{v_p(L)+1}$, at most one term $z+i u$ with $0 \leq i \leq k$ is divisible by $p^{v_p(L)+1}$. Indeed, if two terms $z+i u$ and $z+j u$ with $0 \leq i<j \leq k$ were both divisible by $p^{v_p(L)+1}$, then $p^{v_p(L)+1} \mid(j-i) u$. Since $p \nmid u$, this would imply $p^{v_p(L)+1} \mid j-i$, contradicting $0<j-i \leq k<p^{v_p(L)+1}$. Thus $w_{k-1} \leq v_p(L)$.

For $1 \leq j \leq v_p(L)$, we have $p^j \leq k$. Since $p \nmid u$, the first $p^j$ terms $z, z+u, \ldots, z+\left(p^j-1\right) u$ form a complete residue system modulo $p^j$. Hence at least one term is divisible by $p^j$, so $w_k \geq j$. If $w_{k-1}<j \leq v_p(L)$, then only $w_k$ is at least $j$, and therefore $N_j(z)=1$. Consequently,

\begin{equation}\label{eq:local-form}F_p(n)= \begin{cases}f\left(p^{v_p(d)+v_p(B+A n)}\right)^k, & 0 \leq v_p(B+A n)<v_p(c), \\ f\left(p^{v_p(d)+v_p(c)}\right)^k \prod\limits_{j=1}^{v_p(L)} R_{p, j}^{N_j\left(\frac{B+A n}{p^{v_p(c)}}\right)-1}, & v_p(B+A n) \geq v_p(c) .\end{cases} \end{equation}

Next, we determine the least positive period of $F_p(n)$. Recall that
$$
M_p=\max \left(\left\{j: v_p(k+1)<j \leq v_p(L), R_{p, j} \neq 1\right\} \cup\{0\}\right) .
$$

We distinguish two cases according to the value of $M_p$.

{\bf Case 1.} $M_p>0$. Then $v_p(k+1)<M_p \leq v_p(L), R_{p, M_p} \neq 1$, and $R_{p, j}=1$ for $M_p<j \leq v_p(L)$. Thus the factors with $j>M_p$ in the product in \eqref{eq:local-form} are all equal to 1 and can be omitted. We show that $p^{v_p(c)+M_p}$ is a period.

If $0 \leq v_p(B+A n)<v_p(c)$, then $p^{v_p(B+A n)} \| B+A\left(n+p^{v_p(c)+M_p}\right)$, so

$$
F_p\left(n+p^{v_p(c)+M_p}\right)=F_p(n) .
$$

If $v_p(B+A n) \geq v_p(c)$, write $B+A n=p^{v_p(c)} z$. Then

$$
B+A\left(n+p^{v_p(c)+M_p}\right)=p^{v_p(c)}\left(z+A p^{M_p}\right) .
$$

For $1 \leq j \leq M_p$, the function $N_j$ is $p^j$-periodic and $p^j \mid A p^{M_p}$, so $N_j\left(z+A p^{M_p}\right)=N_j(z)$. Hence \eqref{eq:local-form} gives $F_p\left(n+p^{v_p(c)+M_p}\right)=F_p(n)$. Together with the preceding case, this proves that $p^{v_p(c)+M_p}$ is a period of $F_p$.

Put $T=p^{v_p(c)+M_p-1}$. By Lemma~\ref{lem:counting}, $N_{M_p}$ takes two consecutive integer values and has least positive period $p^{M_p}$. Since $p \nmid A$, the integer $A p^{M_p-1}$ is not a period of $N_{M_p}$. Hence there exists a positive integer $z_0$ such that

$$
N_{M_p}\left(z_0+A p^{M_p-1}\right)-N_{M_p}\left(z_0\right)= \pm 1 .
$$
Choose a positive integer $n_0$ such that
$$
B+A n_0 \equiv p^{v_p(c)} z_0 \quad\left(\bmod p^{v_p(c)+M_p}\right) .
$$
Consequently, $\left(B+A n_0\right) / p^{v_p(c)} \equiv z_0\left(\bmod p^{M_p}\right)$. For $1 \leq j<M_p$, we have $p^j \mid A p^{M_p-1}$, so the corresponding counts $N_j$ are unchanged. Thus \eqref{eq:local-form} gives

$$
\frac{F_p\left(n_0+T\right)}{F_p\left(n_0\right)}=R_{p, M_p}^{N_{M_p}\left(z_0+A p^{M_p-1}\right)-N_{M_p}\left(z_0\right)}=R_{p, M_p}^{ \pm 1} \neq 1 .
$$

Therefore $T$ is not a period. Since the least positive period divides $p^{v_p(c)+M_p}$, it equals $p^{v_p(c)+M_p}=p^{e_p}$.

{\bf Case 2.} $M_p=0$. By the definition of $M_p$, we have $R_{p, j}=1$ for $v_p(k+1)<j \leq v_p(L)$. For $1 \leq j \leq \min \left(v_p(L), v_p(k+1)\right)$, Lemma~\ref{lem:counting} gives $N_j(z)=\frac{k+1}{ p^j}$. Hence, whenever $v_p(B+A n) \geq v_p(c)$, \eqref{eq:local-form} gives

$$
F_p(n)=f\left(p^{v_p(d)+v_p(c)}\right)^k \prod\limits_{j=1}^{\min \left(v_p(L), v_p(k+1)\right)} R_{p, j}^{\frac{k+1}{p^j}-1}=K_p .
$$
If $e_p=0$, then $F_p(n) \equiv K_p$. Suppose that $e_p \geq 1$. By the definition of $e_p$,
$$
f\left(p^{v_p(d)+e_p-1}\right)^k \neq K_p, \quad f\left(p^{v_p(d)+j}\right)^k=K_p \quad\left(e_p \leq j<v_p(c)\right) .
$$

Thus $F_p(n)=K_p$ whenever $v_p(B+A n) \geq e_p$. If $v_p(B+A n)<e_p$, this valuation is determined by $B+A n$ modulo $p^{e_p}$. Hence $F_p(n)$ is determined by $B+A n$ modulo $p^{e_p}$. Since $B+A\left(n+p^{e_p}\right) \equiv B+A n\left(\bmod p^{e_p}\right), p^{e_p}$ is a period of $F_p$. Since $p \nmid A$, choose a positive integer $n_0$ such that $B+A n_0 \equiv 0\left(\bmod p^{e_p}\right)$. Then $v_p\left(B+A\left(n_0+p^{e_p-1}\right)\right)=e_p-1$, and hence
$$
F_p\left(n_0\right)=K_p \neq f\left(p^{v_p(d)+e_p-1}\right)^k=F_p\left(n_0+p^{e_p-1}\right) .
$$

Therefore $p^{e_p-1}$ is not a period, and the least positive period of $F_p$ is $p^{e_p}$.
The least periods $p^{e_p}$ for distinct primes are pairwise coprime. Thus \eqref{eq:finite-product} and Lemma~\ref{lem:product-period} show that the least positive period of $G_f$ is
$$
P_f=\prod\limits_{\substack{p \mid c L \\ p \nmid A}} p^{e_p} .
$$

If the product is empty, then $G_f$ is constant and $P_f=1$.
\end{proof}

\begin{proof}[Proof of Corollary~\ref{cor:monotone}] Fix $p \mid c L$ with $p \nmid A$. Strict monotonicity gives $R_{p, j} \neq 1$ for $1 \leq j \leq v_p(L)$. If $v_p(L)>v_p(k+1)$, then $M_p=v_p(L)$ and $e_p=v_p(c)+v_p(L)$. Suppose that $v_p(L) \leq v_p(k+1)$. Then $M_p=0$. If $v_p(c)=0$, the formula for $e_p$ gives $e_p=0$. Now suppose that $v_p(c)>0$. If $f\left(p^m\right)$ is increasing in $m$, then $R_{p, j}>1$ for $1 \leq j \leq v_p(L)$. The exponents in the product defining $K_p$ are nonnegative, so
$$
K_p \geq f\left(p^{v_p(d)+v_p(c)}\right)^k>f\left(p^{v_p(d)+v_p(c)-1}\right)^k .
$$

If $f\left(p^m\right)$ is decreasing in $m$, then $0<R_{p, j}<1$ for these $j$, and both inequalities are reversed. Thus $e_p=v_p(c)$ in either case, and the result follows from Theorem~\ref{thm:main}.
\end{proof}


\begin{thebibliography}{9}

\bibitem{Farhi}
B. Farhi, \textit{Nontrivial lower bounds for the least common
multiple of some finite sequences of integers},
J. Number Theory \textbf{125} (2007), 393--411.

\bibitem{FarhiKane}
B. Farhi and D. Kane, \textit{New results on the least common
multiple of consecutive integers},
Proc. Amer. Math. Soc. \textbf{137} (2009), 1933--1939.

\bibitem{HongQian}
S. Hong and G. Qian, \textit{The least common multiple of consecutive
arithmetic progression terms},
Proc. Edinb. Math. Soc. \textbf{54} (2011), 431--441.

\bibitem{HongYang}
S. Hong and Y. Yang, \textit{On the periodicity of an arithmetical
function}, C. R. Math. Acad. Sci. Paris \textbf{346} (2008), 717--721.

\bibitem{JiJi}
Q.-Z. Ji and C.-G. Ji, \textit{On the periodicity of some Farhi
arithmetical functions},
Proc. Amer. Math. Soc. \textbf{138} (2010), 3025--3035.

\bibitem{QianTanHong}
G. Qian, Q. Tan and S. Hong, \textit{On the periodicity of a class
of arithmetic functions associated with multiplicative functions},
J. Number Theory \textbf{133} (2013), 2005--2020.

\end{thebibliography}
\end{document}